\documentclass[11pt,reqno]{amsart}
\usepackage{latexsym,amscd,amssymb}
\usepackage{times,mathptm}

\newtheorem{thm}{Theorem}

\theoremstyle{remark}
\newtheorem*{rem-nn}{Remark}

\newcommand{\lk}{\operatorname{Lk}}

\begin{document}
\title[Euler characterisitcs]{
A formula for the Euler characteristics of \\
even dimensional triangulated manifolds}
\author[T. Akita]{Toshiyuki Akita}
\address{Department of Mathematics, Hokkaido University,
Sapporo, 060-0810 Japan}
\email{akita@math.sci.hokudai.ac.jp}

\thanks{Partially
supported by the Grant-in-Aid for Scientific Research (C)
(No.17560054) from the Japan Society for Promotion of Sciences.}

\addtolength{\parskip}{2mm}
\setlength{\parindent}{0.0mm}

\maketitle

A finite simplicial complex $K$ is called an
{\em Eulerian manifold} (or a {\em semi-Eulerian complex}
in the literature) if all of maximal faces have
the same dimension and, for every nonempty face
$\sigma\in K$,
\[
\chi(\lk\sigma)=\chi(S^{\dim K-\dim\sigma-1})
\]
holds, where $\lk\sigma$ is the link of $\sigma$ in $K$ and
$S^n$ is the $n$-dimensional sphere.
Note that $K$ is not necessary connected.
Any triangulation of a closed manifold is an Eulerian manifold.
More generally, a triangulation of a homology manifold
without boundary provides an Eulerian manifold.
The purpose of this short note is to prove the
following alternative formula for the Euler characteristics
of even dimensional Eulerian manifolds.
\begin{thm}\label{main}
Let $K$ be a $2m$-dimensional Eulerian manifold.
Then
\begin{equation}\label{eq-main}
\chi(K)=\sum_{i=0}^{2m}\left(-\frac{1}{2}\right)^i f_i(K)
\end{equation}
holds, where $f_i(K)$ is the number of $i$-simplices of $K$.
\end{thm}
A finite simplicial complex $L$ is called a {\em flag complex} if
every collection of vertices of $L$ which are pairwise
adjacent spans a simplex of $L$.
The formula \eqref{eq-main} was proved in \cite{akita}
under the additional assumptions that $K$ is a PL-triangulation
of a closed $2m$-manifold and is a flag complex.
M. W. Davis pointed out that the formula
\eqref{eq-main} follows from a result in \cite{davis},
% \cite[Remark (ii)]{davis},
provided $K$ is a flag complex
(see {\em Note added in proof} in \cite{akita}).
Both results follow from the considerations of the
Euler characteristics of Coxeter groups.
In this note, we deduce the formula \eqref{eq-main}
from the generalized Dehn-Sommerville equations
proved by Klee \cite{klee}.

Let $K$ be a finite $(d-1)$-dimensional
simplicial complex and $f_i=f_i(K)$
the number of $i$-simplices of $K$ as before.
The $d$-tuple $(f_0,f_1,\dots,f_{d-1})$ is called the
{\em $f$-vector} of $K$.
The {\em $f$-polynomial} $f_K(t)$ of $K$ is defined by
\[
f_K(t)=t^d+f_0t^{d-1}+\cdots+f_{d-2}t+f_{d-1}.
\]
Define the {\em $h$-polynomial} $h_K(t)$ of $K$,
\[
h_K(t)=h_0t^d+h_1t^{d-1}+\cdots+h_{d-1}t+h_d,
\]
by the rule $h_K(t)=f_K(t-1)$.
The $(d+1)$-tuple $(h_0,h_1,\dots,h_d)$ is called the
{\em $h$-vector} of $K$.
%See \cite{ed} for the details of the $f$-polynomial
%and the $h$-polynomial.
The $h$-vector of $K$ satisfies the generalized
Dehn-Sommerville equations, as stated below
in Theorem \ref{DS}.
\begin{thm}[\cite{klee}]\label{DS}
Let $K$ be a $(d-1)$-dimensional Eulerian manifold.
Then
%\begin{equation}\label{eq-DS}
\[
h_{d-i}-h_i=(-1)^i\binom{d}{i}(\chi(K)-\chi(S^{d-1}))
\]
%\end{equation}
holds for all $i$ $(0\leq i\leq d)$.
\end{thm}
\begin{rem-nn}
Klee stated the generalized Dehn-Sommerville equations
in terms of the $f$-vector rather than the $h$-vector.
The formulae quoted in Theorem \ref{DS} are equivalent to
those in \cite{klee} and can be found in \cite{ed}.
Theorem \ref{DS} was also proved in \cite{panov}
by a quite different method,
provided that $K$ is a triangulation of a closed manifold.
\end{rem-nn}
Now we prove Theorem \ref{main}. We have
\[
h_K(-1)=\sum_{i=0}^{2m+1}(-1)^{2m+1-i}h_i
=\sum_{i=0}^{m} (-1)^i (h_{2m+1-i}-h_i).
\]
Now Theorem \ref{DS} asserts that
\[
h_{2m+1-i}-h_i=(-1)^i\binom{2m+1}{i}(\chi(K)-2).
\]
%for $0\leq i\leq 2m+1$.
Hence we obtain
\begin{equation}\label{h-poly}
h_K(-1)=(\chi(K)-2)\sum_{i=0}^m\binom{2m+1}{i}
=2^{2m}(\chi(K)-2).
\end{equation}
On the other hand, we have
\begin{equation}\label{f-poly}
f_K(-2)=(-2)^{2m+1}+\sum_{i=0}^{2m}(-2)^{2m-i}f_i
=2^{2m}\left( -2+\sum_{i=0}^{2m}\left(-\frac{1}{2}\right)^if_i
\right).
\end{equation}
Since $h_K(-1)=f_K(-2)$ by the definition of the $h$-polynomial
$h_K(t)$,
Theorem \ref{main} follows from \eqref{h-poly} and \eqref{f-poly}.

\bibliographystyle{amsplain}

\providecommand{\bysame}{
\leavevmode\hbox to3em{\hrulefill}\thinspace}

\end{document}